\documentclass{amsart}

\usepackage[T1]{fontenc}
\usepackage[utf8]{inputenc}
 \usepackage[main=russian,english]{babel} 

\usepackage{epsfig}
\usepackage{amssymb}
\usepackage{amsmath}
\usepackage{graphicx}

\usepackage{hyperref}
\usepackage{shadow}
\usepackage{boxedminipage}
\usepackage{fancybox}
\usepackage{color}
\usepackage[small,nohug,heads=LaTeX]{diagrams}
\diagramstyle[labelstyle=\scriptstyle]

\title[`Algebra + Computer Science']{A new course `Algebra + Computer Science':\\ What should be its outcomes\\ and where it should start}

\author{Alexandre Borovik}
\thanks{\copyright\ 2022 Alexandre Borovik and Vladimir Kondratiev. Submitted for publication.}
\date{November  24,  2022}
\address{}
\email{alexandre$\gg\! at\! \ll$borovik.net}
\author{Vladimir Kondratiev}
\email{kondratjew239@gmail.com}

\newcommand{\bea}{\begin{eqnarray*}}
\newcommand{\eea}{\end{eqnarray*}}
\newcommand{\bt}{\begin{theorem}}
\newcommand{\et}{\end{theorem}}
\newcommand{\bi}{\begin{itemize}}
\newcommand{\ei}{\end{itemize}}
\newcommand{\bq}{\begin{quote}}
\newcommand{\eq}{\end{quote}}
\newcommand{\bd}{\begin{description}}
\newcommand{\ed}{\end{description}}
\newcommand{\ben}{\begin{enumerate}}
\newcommand{\een}{\end{enumerate}}

\newcommand{\bm}{\setlength{\columnsep}{40pt}\begin{multicols}{2}}
\newcommand{\xbm}{\end{multicols}}

\def\psl{{\rm{PSL}}}

\newcounter{mct}

\newcommand{\Hr}{%
\bigskip
\begin{center}
\begin{minipage}{12.5cm}
  \hrulefill
\end{minipage}
\end{center}
}

\newcommand{\hamburger}{\includegraphics[height=3ex]{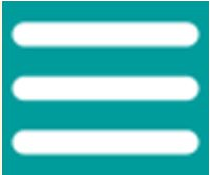}}
\begin{document}

\maketitle
\bq
\textsc{Abstract} The words ``Programming is the second literacy'' were coined more than 40 years ago \cite{Ershov1981}, but never came to life. The paper develops and details  that old slogan by proposing that the mainstream mathematics education in schools should merge with education in  computer science / programming. Of course, this means a deep structural reform of school mathematics education. We are not talking about adapting the 20th century mathematics to the 21st century---as it outlined in \cite{Borovik-Kocsis-Kondratiev2022,Khalin-Vavilov-Yurkov2022}, we mean the 21st century mathematics education for the 21st century mathematics. To the best of our knowledge, this paper is perhaps the first known attempt to start a proper feasibility study for  this reform. The scope of the paper does not allow us to touch the delicate socio-political (and financial) sides of the reform, we are looking only at general curricular and  didactic aspects and  possible directions of the reform. In particular,  we indicate approaches to development of a Domain Specifiic Language (DSL) as a basis for all programming aspects of a new course.
\eq

\section{Executive Summary}

This text discusses some aspects of curriculum for a  proposed replacement of the traditional school algebra course by a new course which fully integrates algebra with  informatics / programming. Approaching this issue from the basic principles of project management, we focus on two points from which development of any project starts:
\ben
\item What are our aims? What do we wish to achieve?

\item What is our starting position?
\een

In the context of curriculum development, this becomes

\bd
\item[Question 1] What are learning outcomes of the new course? What school graduates should be able to do on completion of the course?
\item[Question 2] How the new course will relate to the primary school arithmetic?
\ed

In case of projects directed at work with children for considerable periods of time (in the case of algebra curriculum, this is perhaps 7 or even more years), point (2) and Question 2 attain larger importance than, say, in most industrial projects: the child has to be in the focus of the project, and the child grows, and all the proposed activities have to grow with him.

Both points need to be clearly described and illustrated by sample problems of the kind Learner is supposed to be able to solve (on their own, without guidance from teachers) at the start and at the end of their  study. Only then we can start writing a curriculum, working  in the direction from the answer to Question 1 to answer to Question 2, at every step including only the material which is definitely needed for  further progression  to a higher level.

\subsection{} \textbf{Regarding Question 1, we formulate the following recommendations.}

\subsubsection{}  The most important issue of the 21st Century is the relationship between people and computers.

Therefore we expect that on completing the course, Learner should
\bi
\item have understanding of opportunities and dangers coming from saturation of every aspect of the economy, politics, warfare, media, entertainment, everyday life, by computers and computer systems which make decision on behalf of people;
\item in navigating this dangerous world, understand what questions could and should be asked---and to whom;
\item if necessary, be prepared to learn, and use, more technical and professional aspects of computer science and computer programming.
\ei
This could be compared with mathematical skills expected  from secondary school graduates in Europe in the first half of the 20th century: they had to master enough of algebra and trigonometry for their subsequent training, if needed, as artillery officers or air pilots---after all, it was the era of mass conscription armies.

\subsubsection{} There is a simple and challenging criterion of Learner's the mastery of the required skills.
\bi
\item From a relatively early  stage of the course, Learner's answer to an arithmetic or algebraic problem (including mathematical problems described as `real life' problems) should be an executable computer code developed, almost in its entirety, by Learner---without use of standard packages such as \textsc{Mathematica}---which
\bi
\item solves \emph{all} problems of the same type;
\item helps to check, analyse, and generalise the solution.
\ei
Problems which require a \emph{proof} of some statement are naturally excluded from this requirement---but computer experiments should be welcome in construction of counterexamples and in formulation of conjectures.
\ei

\subsubsection{} The requirements for a computer language to be used in the course are quite demanding.
\bi
\item If more advanced versions of it are used at later stages of learning, it needs to be backward compatible at all stages of school algebra---to allow a systematic revisiting, reusing, and rethinking of earlier learned material, problems solved, and codes written.
\item But its small ``starting''  fragment should provide a simple, safe,  and exiting playground for children at the first steps of learning algebra.
\item The language should include sufficiently rich elements usually found in high level languages and provide efficient preparation for learning professional industry standard languages.
\item  All numbers have to be exclusively symbolic, including fractions and surds; if the result in some intermediate calculation  is $\frac{1-\sqrt{3}}{2}$, it stays that way and is used in subsequent calculation.  Of course, they have to include constants $e$  and $\pi$. No rounding, no floating point decimals.\footnote{And the answer to the question ``Find the reminder upon division of polynomial $ x^{2023} + 1$ upon $x^2 - 4$'' should be $2^{2022}\cdot x +1$, without any attempt to print $2^{2022}$ in decimals.}
\ei

This will require development of a bespoke \emph{Domain Specific Language}\footnote{From Wikipedia: ``A Domain-Specific Language (DSL) is a computer language specialized to a particular application domain. This is in contrast to a general-purpose language (GPL), which is broadly applicable across domains.'' Famous examples include  \textsf{html} and \LaTeX. Some books: \cite{Boesma2020,Fowler2010,Voelter2010}. We are using Wikipidea as a primary reference because neither us, nor our expected readers are computer scientists or computer programmers.} (DSL) for use in the course. Most likely, existing general purpose languages are not suitable for this role.

\begin{center}
\textbf{This texts represents first early steps in development of specifications for DSL.}
\end{center}

\subsubsection{} The content of the algebra course should change and reflect the demands of computer programming, and include, for example, Boolean algebra, elements of number theory, modular arithmetic, and finite fields. Shaping the algebraic curriculum should go in step with development of the DSL.

\subsubsection{} Other parts of mathematics, for example, geometry, mechanics, and statistics need a separate discussion and are not touched here. Also, interactions with physics deserve most serious attention.

\subsection{}   \textbf{On Question 2, our principal suggestions are the following:}

\bi
\item Introduction of ``computational'' thinking  in algebra should be prepared by development of algorithmic thinking in arithmetic.

\item ``Questions method'' for word problems in arithmetic is suggested as an efficient tool for mastering an elementary level of  algorithmic thinking---and this has done in arithmetic before moving to algebra.

\item   Use of a  simple computer language (which still needs to be developed) is essential from the first day of learning algebra, in particular, for formulation of  algorithms developed by the arithmetic `questions method' in a compact form suitable for conversion into computer code.

\item Typed variables and type inference are essential for a computer language that will be used in the course, and their introduction needs to be prepared by mastering `named numbers' in arithmetic. It needs to be clarified that
    \bq
   \emph{named numbers are \textbf{not} numbers replaced by names and symbols,}
    \eq
they are numbers of units of various kind: 10 apples are not the same as 10 people and not the same as 10 kilometers, and you do \emph{not} add 10 apples and 10 people.

\item However, there is a conceptual gap between
\bi
\item formulation of an algorithm and its implementation in a code,  and
\item representation of the answer as a closed algebraic formula.
\ei
Therefore  Learner will still need mastering algebraic manipulations with formulae in algebra---which is conceptually parallel to use of evaluation by call by name in computer programming---as opposed to evaluation by call by value used in a  direct implementation of an a arithmetic algorithm.
\item A lot of attention needs to be given to \emph{arithmetic of named numbers}. In particular,  Learner should learn that they can freely introduce in their solutions intermediate parameters, and their values, if the names of these parameters (correspondingly, types of variables) are alien, appear neither in the data (conditions) of the problem, no in the supposed answer. These additional parameters will disappear in the course of solution.
\item This approach inevitably requires  Learner reflect on his/her work, analyse it, apply, in effect, some meta-thinking (that is, thinking about thinking)---and can be achieved only if  Learner is psychologically comfortable and feels being in full control of the problem, its solution, and the computer.
\ei
\medskip

\begin{center}
\shadowbox{\parbox{10cm}{%
Achieving all that is a challenging task. We would not advise to undertake it lightly.
}}
\end{center}

\medskip

\subsection{Warning about the scale and cost of the project} \label{sec:scale} \

The first author has experience of working on large and complex research projects, see  for example \cite{ABC}, a book of 550+ pages that he co-authored and which contained a proof  of a single theorem---the proof was in making for about 15 years and used crucially important contributions from an informal team of 10 people. Closer to the theme `algebra + computer science for school' is his work on computational symbolic logic, where his joint paper with \c{S}\"{u}kr\"{u} Yal\c{c}\i nkaya  \cite{BY2018} was developed, step by step, over at least 10 years, and progress  critically depended  on systematic computer experiments.

Basing his assessment on this experience, he is quite confident that  development of a software system supporting the proposed course and \emph{used on the scale of a nation} is realistic, but could easily cost up to  a few million dollars, demand work of a large interdisciplinary team of experts,  and, with necessary test runs in schools, could take up to ten years. This estimate does not include the cost of re-education of the whole army of teachers, and time required for that.

We urge the reader to remember the words of H. L. Mencken:
\bq
\emph{For every complex problem there is an answer that is clear, simple, and wrong.}
\eq

We warn:

\medskip

\begin{center}
\shadowbox{\parbox{10cm}{%
Beware of snake oil merchants peddling cheap and easy recipes for revitalising school mathematics by means of `digital technologies'!
}}
\end{center}

\medskip

We hope that our work could be used as an antidote against their promises.

\medskip

\begin{center}
  ***
\end{center}

\medskip

The rest of the paper contains a more detailed discussion of Question 2:
\medskip

\begin{center}
\textbf{How the proposed  course will relate to the primary school arithmetic?}
\end{center}

\medskip

\section{But what is arithmetic?}

School arithmetic contains two major parts.
\bi
\item Formal written methods for arithmetic operations on decimals and fractions \cite{Gardiner2016}, that is, manual computation based on certain recursive algorithms.
\item Solving word problems.
\ei
As explained in \cite{Borovik2017-Simplicity}, ``word problems” of arithmetic
involve identification of mathematical structures and relations of the real world and
mapping them onto better formalised structures and relations of arithmetic, or, in
Igor Arnold's formulation  from 1946, when the words `structure' and `relation' were not yet \emph{en vogue} \cite{Arnold1946},
\bq
 \dots\  \emph{teaching arithmetic involves, as a key component, the
development of an ability to negotiate situations whose concrete natures represent very different relations between magnitudes and quantities.}
\eq

Even more important is Igor Arnold’s characterisation of arithmetic:
\bq
\emph{The difference between the ``arithmetic” approach to solving problems and the algebraic
one is, primarily the need to make a concrete and sensible interpretation of all the values
which are used and/or which appear at any stage of the discourse.}
\eq

Therefore the phase transition from algebra to arithmetic is an important turning point in learning mathematics.

\section{The bridge from arithmetic to algebra}
\label{sec:bridge}

This is our principal thesis:
\medskip

\begin{center}
\shadowbox{\parbox{10cm}{%
The merger of school mathematics with informatics should start simultaneously with the phase transition of mathematical learning from arithmetic to algebra, and should be rooted in those aspects of arithmetic which actually belong to computer science but are not usually recognised as such.
}}
\end{center}
\medskip\noindent
For example, arithmetic already contains
\bd
\item[Abstraction] The concept of \emph{number} is already a huge abstractio.

\item[Algorithms] First of all, long division, etc. are algorithms. But they are given to children as something God-sent, just as rules to follow. However, what will be shown here, the  `questions method' of yesteryear arithmetic (now almost universally forgotten (one of its rare  clear exposition is \cite{Beloshistaya2007}) allows children to develop their own algorithms which solve many types of arithmetic problems.

\item[Recursion] The notorious long division is a recursive algorithm, as well as long multiplication, as well as addition and subtraction of  decimals. Euclid's algorithm for finding the greatest common divisor of two integers is also recursive.

\item[Typed variables] `Named numbers' of arithmetic  perfectly map to typed variables of computer programming.

\item[Reification] This is a high level concept from computer science, but its toy version appears in arithmetic as introduction of intermediate parameters (`helpful numbers' is a  better term for kids' use), see Section \ref{sec:auxiliary}. Or as 'fantasy  units of measurement' -- see \cite[Sections 3.1 and 3.3]{Borovik-Shadows} and \cite[Solutions and Notes to Problems 89 and 90]{Borovik-Gardiner2019} (in the latter they are called `hidden parameters') for discussion and historic examples.

     In computer science, \emph{reification} is a widely used concept. This is how \textsc{Wikipedia} defines reification:
\bq
\emph{Reification} is the process by which an abstract idea about a computer program is turned into an explicit data model or other object created in a programming language. A computable/addressable object -- a resource -- is created in a system as a proxy for a non computable/addressable object.\footnote{\url{https://en.wikipedia.org/wiki/Reification_(computer_science)}.}
\eq
\noindent
In school algebra, the notorious $x$ as a notation for the unknown is a typical example of reification: it is a proxy for a number  not given to us, but a proxy  nevertheless allowing algebraic manipulations with it.

\ed

The rest of this text concentrates on this issue: the ``questions method'' and algorithms in transition from arithmetic to algebra. Its principal point  is

\smallskip

\begin{center}
\shadowbox{\parbox{10cm}{%
Fusion of algebra and computer science should be prepared by development of algorithmic thinking in arithmetic.
}}
\end{center}

We shall explain, that, in school arithmetic of yesteryear, there was an  approach to an efficient development of algorithmic thinking: the so-called ``questions method'' for solving word problems\footnote{The questions method is discussed in some details in \cite[Sections 3 and 4]{Borovik2017-Simplicity}.}. We further claim that it could serve as a bridge between arithmetic and algebra, especially if algebra is merged with computer science  / coding.

We shall also try to demonstrate that the questions method allows us to reach a number of didactic  targets:

\bi
\item Emphasis on detection and analysis of structures and relations of the real world and their representation in mathematical terms.
\item Systematic analysis of data given in a problem.
\item Correspondence between named numbers in algebra and typed  variables in code writing.
\item An answer to an arithmetic or algebraic problem being an executable computer code which
\bi
\item solves \emph{all} problems of the same type;
\item helps to check and analyse the solution;
\item uses evaluation by call by name and by call by value.
\ei
\ei
If these targets are achieved, we will be able to claim that

\smallskip

\begin{center}
\shadowbox{\parbox{10cm}{%
The child is in control of the problem, of its solution, and of the computer \cite{Borovik2017-In-control}.
}}
\end{center}

And we have to achieve that---otherwise, no further progress is possible.

We do not discuss here why a fusion of algebra with computer science could be desirable in school education -- this was  discussed in \cite{Borovik-Kocsis-Kondratiev2022} and \cite{Khalin-Vavilov-Yurkov2022}.

\section{A sample problem from arithmetic: Bowl of cherries}\

Let us consider a  simple problem.\footnote{The problem, and the picture illustrating it, was picked from an online  talk by  Olga Moskalenko  «Проблема систематических ошибок в решении текстовых задач по результатам мониторинга на платформе «Учи.ру»» on Научный семинар по методике преподавания математики, Moscow State University, 22 September 2022. Interestingly, in a test only $39\%$ of Year 6 students gave the correct answer, $6$ minutes.} Of course, at school children should start with much simpler, 1 or 2 step problems; problems that the one below should come at a later stage. We use this example because it provides  a wider overview of the arithmetic / algebra boundary.

\Hr
\begin{center}
\textsc{Bowl of cherries}\\
Alice and Bob pick cherries.

\vspace{-14ex}

\parbox[c]{1in}{\vspace{1.5in} Bob: I can fill this bowl in $8$ minutes.}
\parbox[t]{2in}{\begin{center}
\includegraphics[width=1.5in]{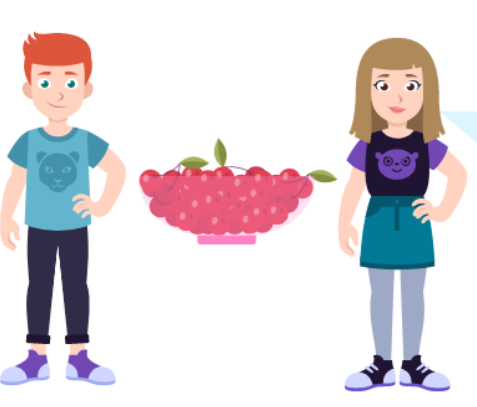}
\end{center}}
\parbox[c]{1in}{\vspace{1.5in}Alice: I can fill this bowl in $24$ minutes.}

\hfill\hfill \hspace{1em}Bob \hspace{2em}\hfill Alice\hspace{2em} \hfill\hfill\emph{}

Working together, in what time they will fill the bowl?
\end{center}
\Hr
\bigskip
\noindent
Let us approach this problem step by step in compliance with the methodology that we have just formulated.

\section{Analysing data} \label{sec:Analysing} \

\subsection{Parsing of the problem} \label{sec:parsing}
Let us now  to ask the all-important question:
\bq
\textbf{What is the data and what is the question?}
\eq
\medskip

\noindent
Indeed, this question is well justified because the  correct parsing of the problem is the following:

\medskip

\begin{center}
\fbox{\parbox{10cm}{
\bd
\item[{Data}] Bob fills the bowl in $8\; \mbox{\texttt{min}}$, Alice in  $24\; \mbox{\texttt{min}}$.\\
$\mbox{ }\quad$ $\mbox{ }$   Now they work together.

\item[{Question}] In what time they will fill the bowl?
\ed
}}
\end{center}

\medskip

We shall call the question in the problem that we separated from the data the \emph{target question}.

\subsection{Questions method and step questions}  And now we shall give an example of a questions method. In the 1960s in Russia, students were supposed to master it by the age of 11 or 12 and be able to produce, for a word problem, a short  sequence of simpler \emph{step questions}, which resulted in answering the target question. This could also be described as a formulation of a multi-step solution to a problem as a composition of solutions to much simpler one step problems.

One of many possible step question sequences for the \textsc{Bowl of cherries} problem is shown in Section \ref{sec:Questions}. We shall try to explain now how it could be obtained.

\subsection{Generic guiding questions}

Teaching children to formulate step questions in a  questions method  focuses on development of each child’s
ability to start his/her ``questions” attempt at a solution by asking himself or herself appropriate self-directed \emph{guiding}  questions (they were called \emph{auxiliary questions} in the Russian pedagogical literature, but in England, the words ``\emph{an auxiliary question}” are
loaded with expectation that the question is asked by a teacher to help a struggling
pupil).

Generic questions (or meta-questions) recommended by the questions method could be data-centered \cite{Alekseeva-Istchenko2019,Bazhenova-Odoevtseva2012,Beloshistaya2007,YaKlass2022}, for example,
\bi
\item What questions can be asked about these data?\footnote{In the 1960s in Russia these questions could be asked even in a harsher form: a teacher wrote a problem on the blackboard, gave a  minute to students to think, then wiped out the target question and asked: what questions could be asked about these data? What could be learnt from these data?}
\item What can be learnt from these data?
\item And a more specific question: how to express mathematically ``working together''?
\ei
Beloshistaya \cite[p. 303]{Beloshistaya2007} calls this approach \emph{synthetic}. It is interesting to compare Beloshistaya's analysis with Gerofsky \cite{Gerofsky1996}.

There are also more advanced  guiding questions. Those found in the modern literature  appear to come directly from a project management manual \cite{YaKlass2022}.
\bq
\emph{Analyse the task working backwards from the target back to the current situation}.
\eq
Indeed look for yourself at questions of the synthetic approach:
\bi
\item What additional quantities it would be useful to know  for finding the answer? This question could lead to introduction of a `helpful number' (an intermediate parameter), see  Section~\ref{sec:auxiliary}.
\item What do you have to do for finding these additional quantities?
\item How can they be deduced from the data given?
\ei

\subsection{The crucial point: To whom these questions were addressed?}

In the old tradition of the ``questions method'',
\bi
\item At the first stages of learning the questions method,  guiding questions  were addressed by a teacher or a textbook to \textbf{students}, aged 11--12.
\item At later stages, \textbf{students} were expected to formulate there guiding questions \textbf{themselves}, thus developing  \textbf{their own step questions} directed at solving a specific problem.
\item Put together, \textbf{student's} step questions formed in fact  an \textbf{algorithm} (although this word was not mentioned): any problem of the same type, but with different data, could be quickly answered by following the same steps.
\item This algorithm had to be developed by the student.
\ei

\subsection{Introduction of an intermediate parameter (or a ``helpful number'')} \label{sec:auxiliary}
Answering a guiding question
\bq
How to express mathematically ``working together''?
\eq
after some thought or discussion it could be discovered that  ``to work together'' means, in mathematical terms, that
\bq
productivities, or speeds of picking cherries, \emph{add up}.
\eq
This should immediately trigger further questions:

\bi
\item How to measure the speed of picking?

Of course, in \emph{\texttt{cherries per minute}}!

\item How many cherries we have in the problem?
\ei
And the last question leads to a \emph{critically important step}:

\medskip
\fbox{\parbox{10cm}{
\
\noindent
\textbf{A helpful number:} Let us assume that the bowl contains $72$ \texttt{cherries}.
}}

\medskip

This choice, $72$, appears to be  arbitrary and made just for the convenience of division by $8$ and by $24$, which allows to avoid fractions in further calculation. We shall soon see that this is well justified. Also, as it has already been mentioned in Section~\ref{sec:bridge}, it is an example of \emph{reification}, one of concepts of computer science already present, in disguised or hidden form, in arithmetic. Later, in Section \ref{sec:elimination}, it will be explained why this step is safe.

\section{Solution} \

\subsection{The sequence of step questions} \label{sec:Questions} The following solution is a sequence of questions and  answers formulated on the basis of analysis of the data and in response to generic  guiding questions described in Section~\ref{sec:Analysing}.

\medskip

\fbox{\parbox{12cm}{
\noindent
Question 1. What is Alice's picking speeds?
\[
3 \frac{\mbox{ \texttt{cherries}}}{\mbox{\texttt{min}}} = \frac{72 \mbox{ \texttt{cherries}} }{24\; \mbox{\texttt{min}}}
\]
Question 2. What are is Bob's picking speed?
\[
9 \frac{\mbox{ \texttt{cherries}}}{\mbox{\texttt{min}}} = \quad \frac{72 \mbox{ \texttt{cherries}} }{8\; \mbox{\texttt{min}}}
\]
Question 3. What is their picking speed working together?
\[
12\frac{\mbox{ \texttt{cherries}}}{\mbox{\texttt{min}}} = 3 \frac{\mbox{ \texttt{cherries}}}{\mbox{\texttt{min}}} + 9 \frac{\mbox{ \texttt{cherries}}}{\mbox{\texttt{min}}}
\]
Question 4. In what time they will fill the bowl?
\[
6 \;\mbox{\texttt{min}} = 72 \mbox{ \texttt{cherries}} \div 12\frac{\mbox{ \texttt{cherries}}}{\mbox{\texttt{min}}}
\]
}}

\subsection{Elimination of intermediate parameters of independent types} \label{sec:elimination}

But what will happen with the answer to the problem if we change the number of cherries---it was not given to us, it was selected by us just for convenience?
\medskip

Nothing will happen.
\medskip

The data given is of type \texttt{min}, the answer is also of type  \texttt{min}. The extra parameter is of completely independent type \texttt{cherry}. The answer does not depend on the numerical value of the extra parameter, because its change can be treated simply as the change of \emph{unit of measurement}: instead of 72 cherries we can talk about 32 \texttt{pairs},  or 24 \texttt{spoons}, or 12 \texttt{cups} of cherries, etc.---all these different numbers describe the same physical quantity.

This simple observation gives us
\begin{center}
\shadowbox{\parbox{10cm}{%
\textbf{Freedom of choice of intermediate parameters} --- we need only to make sure that their types are independent from types of data.
}}
\end{center}

\subsection{The resulting (pseudo)code}

This is the setup of the problem which uses letters instead of numbers:

\bq
 Bob fills the bowl in \texttt{B}\; \texttt{min}, Alice in  \texttt{A}\; \texttt{min}\\
 Working together, they fill the bowl in T\; \texttt{min}.
\eq

A ``questions method'' solution can be easily rewritten as a pseudocode:

\bigskip

\fbox{\parbox{12cm}{
\noindent
\texttt{input A, B \hfill \% Alice's and Bob's times}\\
\texttt{input C \hfill \% number of cherries}

\texttt{U := C:A; \; V :=  C:B;    \hfill \% Alice's and Bob's speeds}

\texttt{W := U + V; \hfill \%  Their picking speed working together}

\texttt{T := C:W \hfill \%  The time in which they fill the bowl}

\noindent
\texttt{return T};
}}

\bigskip

But how to simplify it to a one line pseudocode: \\

\fbox{\parbox{12cm}{
\noindent
\texttt{input A, B;} \quad  \texttt{T:= (A*B):(A+B)};  \quad \texttt{return T}
}}

\medskip
\noindent
that is,
 \[
 T := \frac{AB}{A+B}?
 \]
This is a serious issue.\footnote{In the discussion of this problem at the Association of Teachers of Mathematics (UK) meeting in October 2022, one of my colleagues quite rightly asked ``I am wondering how a student might become aware that $(AB)/(A+B)$ might be a solution to the problem?'' And someone else suggested ``I think of it as $1/(1/A+1/B)$ rather than it's simplified form''.}
 In the early algebra, writing code is \textbf{much} easier than developing an algorithm.  Indeed the answer
 \[
 T := \frac{AB}{A+B}
 \]
---the shortest expression of the algorithm ---is non-trivial.

\subsection{Solution in symbolic variables} \label{sec:solution-symbolic}

Compression of the code obtained by the questions method into a compact formula could be done either manually, or maybe even automatically; for the latter case, the code has to be embedded into a symbolic algebra system which automatically simplifies algebraic expressions. We think it would be desirable to achieve a didactically efficient balance of the two approaches.
\bigskip

\fbox{\parbox{12cm}{
\noindent
Setup: $A$: Alice's speed, $B$: Bob's Speed, $C$: number of cherries.

\medskip
\noindent
Question 1. What are Alice's and Bob's picking speeds?
\[
\frac{C}{A} \quad\mbox{ and }\quad \frac{C}{B}
\]
Question 2. What is their picking speed working together?
\[
V := \frac{C}{A} + \frac{C}{B}
\]
Question 3. In what time they will fill the bowl?
\[
T := \frac{C}{V} == \frac{C}{\frac{C}{A} + \frac{C}{B}} == \frac{AB}{A+B }
\]
}}

\medskip
Notice that this is \emph{``call by name''} evaluation.

Also notice that variable $C$ disappeared from the answer---as it was inspected, because it was of type  independent from the type of data.\footnote{In an earlier discussion of this text, there was an interesting question about this problem: ``Who ate all the cherries?''. Well, we  can add to the setup of the problem someone called Kevin, who eats cherries at the rate of $1$ bowl per $K$ minutes. The code could be easily adjusted producing the answer that the bowl will be filled in time
\[
\frac{ABK}{AK+BK-AB},
\]
and also providing a warning that if
\[
K \leqslant \frac{AB}{A+B},
\]
the bowl will never be filled. However, some colleagues suggested, quite rightly, that the equivalent inequality which characterises the failure to fill the bowl
\[
\frac{1}{K} \geqslant \frac{1}{A} + \frac{1}{B}
\]
was easier to interpret in real life terms.
}

\subsection{A completely different solution}
Of course, the questions method can produce other solutions as well. Here is one of them.

\noindent
Question 1. In how many times Bob is faster than Alice?
\[
\frac{24\; \mbox{\texttt{min}}}{8 \;\mbox{\texttt{min}}} = 3   \quad\mbox {or }\quad  \frac{A}{B}.
\]
\noindent
Question 2. In how many times the two of them working together, are faster than Alice working alone?
\[
3+1 = 4  \quad\mbox {or }\quad \frac{A}{B} + 1.
\]
Question 3. In what time they will fill the bowl?
\[
24\;\mbox{\texttt{min}} : 4 = 6 \;{\texttt{min}} \quad\mbox {or }\quad \frac{A}{\frac{A}{B} + 1} = \frac{AB}{A+B}.
\]

\section{User's Story}

We hope that this ``questions method” solution can be naturally converted into a
computer code with the help of a child friendly GUI interface.  What follows is a brief discussion of a desirable functionality of this interface, a ``User's Story'', in terminology of software development, written in accordance with recommendations of \emph{User experience design}\footnote{User experience design, \url{https://en.m.wikipedia.org/wiki/User_experience_design}.}.

Assuming that a student (we call him/her Learner) types everything in some computing device (smartphone, tablet, laptop) with the course software installed, as soon a number is entered, the GUI starts a dialogue with Learner asking whether this number is  a datum, or a helpful number, or the answer to the problem and what is its type. Also, in this dialogue a number of various ambiguities are resolved---we omit these details here. In simple problems, this dialogue is likely to be very short---we shall explain in the next our paper.

The GUI assigns internal variables to all numbers on the screen, and converts numbers into cells whose values could be changed by the code or edited by Learner The values of data and helpful numbers remain editable by Learner see Figure \ref{fig:screen-1}. It is an executable code represented in a GUI.   The `Hamburger' icon  \hamburger\  is a pop-up menu with many additional functions.

If Learner chooses button {\color{green} \Ovalbox{Run}}  in the menu and presses it, the cells become alive and the GUI looks as in Figure~\ref{fig:screen-2}.

The Learner is invited to experiment with changing data and auxiliary parameters; it the helpful number is changed from 72 to 48, and the button {\color{green} \Ovalbox{Run}} is pressed, the Learner gets another screen,  Figure~\ref{fig:screen-3}---with exactly the same answer!

Even more important: it should be possible to substitute letters in the cells (or in just one cell) for data and helpful numbers, thus automatically getting a symbolic solution of the kind offered in Section \ref{sec:solution-symbolic}. See Figure \ref{fig:screen-4} for what happens if the  time for Alice to fill the bowl is denoted by letter $A$.

The course software should help Learner to move from arithmetic to algebra: a solution of an arithmetic problem can be converted into a code which solves the same problem, but with symbolic inputs.

\begin{figure}
\begin{center}
\footnotesize
\bigskip
\begin{center}
\begin{boxedminipage}[t]{12cm}
\begin{center}
{\color{blue}{\large\textsc{Bowl of cherries}}}
\end{center}
\end{boxedminipage}
\end{center}
\begin{center}
\begin{boxedminipage}[t]{12cm}
\

\textbf{Data:} Bob fills the bowl in {\color{black} \Ovalbox{$8\; \mbox{\texttt{min}}$}}, Alice in {\color{black} \Ovalbox{$24\; \mbox{\texttt{min}}$}}.\\
$\mbox{ }\quad$ $\mbox{ }$   Now they work together.

\textbf{Question:}  In what time they will fill the bowl?
\end{boxedminipage}
\begin{boxedminipage}[t]{12cm}

\textbf{Helpful number: }The bowl contains {\color{black} \Ovalbox{$72\; \mbox{\texttt{cherry}}$}}.

\end{boxedminipage} \\

\begin{boxedminipage}[t]{12cm}
\

\noindent
\textbf{Question 1.} What is Alice's picking speeds?

\textbf{Answer:} \Ovalbox{\color{black}$3 \frac{\mbox{ \texttt{cherry}}}{\mbox{\texttt{min}}}$} = \Ovalbox{\color{black} $72 \mbox{ \texttt{cherry}}$} $\div$ \Ovalbox{\color{black}$24\; \mbox{\texttt{min}}$}

\textbf{Question 2.} What are is Bob's picking speed?

\textbf{Answer:}  \Ovalbox{\color{black}$9 \frac{\mbox{ \texttt{cherry}}}{\mbox{\texttt{min}}}$} = \Ovalbox{\color{black} $72 \mbox{ \texttt{cherry}}$} $\div$ \Ovalbox{\color{black}$8\; \mbox{\texttt{min}}$}

\textbf{Question 3}. What is their picking speed working together?

\textbf{Answer:}  \Ovalbox{\color{black}$12\frac{\mbox{ \texttt{cherry}}}{\mbox{\texttt{min}}}$} = \Ovalbox{\color{black}$3 \frac{\mbox{ \texttt{cherry}}}{\mbox{\texttt{min}}}$} $+$ \Ovalbox{\color{black}$9 \frac{\mbox{ \texttt{cherry}}}{\mbox{\texttt{min}}}$}

\textbf{Question 4.} In what time they will fill the bowl?

\textbf{Answer:}  \Ovalbox{\color{black}$6 \;\mbox{\texttt{min}}$} = \Ovalbox{\color{black}$72 \mbox{ \texttt{cherry}}$} $\div$ \Ovalbox{\color{black}$12\frac{\mbox{ \texttt{cherry}}}{\mbox{\texttt{min}}}$}
\hfill  \hamburger
\end{boxedminipage}
\end{center}
\caption{GUI after entering the solution. The `Hamburger' icon  \hamburger\  is a pop-up menu with very additional functions.}
\label{fig:screen-1}
\normalsize
\end{center}
\end{figure}

\normalsize


\begin{figure}
\begin{center}
\footnotesize
\bigskip
\begin{center}
\begin{boxedminipage}[t]{12cm}
\begin{center}
{\color{blue}{\large\textsc{Bowl of cherries}}}
\end{center}
\end{boxedminipage}
\end{center}
\begin{center}
\begin{boxedminipage}[t]{12cm}
\

\textbf{Data:} Bob fills the bowl in {\color{red} \Ovalbox{$8\; \mbox{\texttt{min}}$}}, Alice in {\color{red} \Ovalbox{$24\; \mbox{\texttt{min}}$}}.\\
$\mbox{ }\quad$ $\mbox{ }$   Now they work together.

\textbf{Question:}  In what time they will fill the bowl?
\end{boxedminipage}
\begin{boxedminipage}[t]{12cm}

\textbf{Helpful number: }The bowl contains {\color{red} \Ovalbox{$72\; \mbox{\texttt{cherry}}$}}.

\end{boxedminipage} \\

\begin{boxedminipage}[t]{12cm}
\

\noindent
\textbf{Question 1.} What is Alice's picking speeds?

\textbf{Answer:} \Ovalbox{\color{blue}$3 \frac{\mbox{ \texttt{cherry}}}{\mbox{\texttt{min}}}$} = \Ovalbox{\color{red} $72 \mbox{ \texttt{cherry}}$} $\div$ \Ovalbox{\color{red}$24\; \mbox{\texttt{min}}$}

\textbf{Question 2.} What are is Bob's picking speed?

\textbf{Answer:}  \Ovalbox{\color{blue}$9 \frac{\mbox{ \texttt{cherry}}}{\mbox{\texttt{min}}}$} = \Ovalbox{\color{red} $72 \mbox{ \texttt{cherry}}$} $\div$ \Ovalbox{\color{red}$8\; \mbox{\texttt{min}}$}

\textbf{Question 3}. What is their picking speed working together?

\textbf{Answer:}  \Ovalbox{\color{blue}$12\frac{\mbox{ \texttt{cherry}}}{\mbox{\texttt{min}}}$} = \Ovalbox{\color{blue}$3 \frac{\mbox{ \texttt{cherry}}}{\mbox{\texttt{min}}}$} $+$ \Ovalbox{\color{blue}$9 \frac{\mbox{ \texttt{cherry}}}{\mbox{\texttt{min}}}$}

\textbf{Question 4.} In what time they will fill the bowl?

\textbf{Answer:}  \Ovalbox{\color{blue}$6 \;\mbox{\texttt{min}}$} = \Ovalbox{\color{blue}$72 \mbox{ \texttt{cherry}}$} $\div$ \Ovalbox{\color{blue}$12\frac{\mbox{ \texttt{cherry}}}{\mbox{\texttt{min}}}$}

\end{boxedminipage} \\\begin{boxedminipage}[t]{12cm}
\

\begin{center}
 The answer to the problem is\\
Working together, Alice and Bob will fill the bowl in {\color{green} \Ovalbox{$6\; \mbox{\texttt{min}}$}}\\

\end{center}

\end{boxedminipage} \\

\begin{boxedminipage}[t]{12cm}
\
You may click on highlighted red {\color{red}\Ovalbox{numbers}} and change them! You will get a different problem of the same type.
\hfill\hamburger

\end{boxedminipage} \\
\end{center}

\caption{The code goes alive. This is something that could be proudly shown to a parent.}
\label{fig:screen-2}
\normalsize
\end{center}

\end{figure}

\begin{figure}
\begin{center}
\footnotesize
\bigskip
\begin{center}
\begin{boxedminipage}[t]{12cm}
\begin{center}
{\color{blue}{\large\textsc{Bowl of cherries}}}
\end{center}
\end{boxedminipage}
\end{center}
\begin{center}
\begin{boxedminipage}[t]{12cm}
\

\textbf{Data:} Bob fills the bowl in {\color{red} \Ovalbox{$8\; \mbox{\texttt{min}}$}}, Alice in {\color{red} \Ovalbox{$24\; \mbox{\texttt{min}}$}}.\\
$\mbox{ }\quad$ $\mbox{ }$   Now they work together.

\textbf{Question:}  In what time they will fill the bowl?
\end{boxedminipage}
\begin{boxedminipage}[t]{12cm}

\textbf{Helpful number: }The bowl contains {\color{red} \Ovalbox{$48\; \mbox{\texttt{cherry}}$}}.

\end{boxedminipage} \\

\begin{boxedminipage}[t]{12cm}
\

\noindent
\textbf{Question 1.} What is Alice's picking speeds?

\textbf{Answer:} \Ovalbox{\color{blue}$2 \frac{\mbox{ \texttt{cherry}}}{\mbox{\texttt{min}}}$} = \Ovalbox{\color{red} $48 \mbox{ \texttt{cherry}}$} $\div$ \Ovalbox{\color{red}$24\; \mbox{\texttt{min}}$}

\textbf{Question 2.} What  is Bob's picking speed?

\textbf{Answer:}  \Ovalbox{\color{blue}$6 \frac{\mbox{ \texttt{cherry}}}{\mbox{\texttt{min}}}$} = \Ovalbox{\color{red} $48 \mbox{ \texttt{cherry}}$} $\div$ \Ovalbox{\color{red}$8\; \mbox{\texttt{min}}$}

\textbf{Question 3}. What is their picking speed working together?

\textbf{Answer:}  \Ovalbox{\color{blue}$8\frac{\mbox{ \texttt{cherry}}}{\mbox{\texttt{min}}}$} = \Ovalbox{\color{blue}$2 \frac{\mbox{ \texttt{cherry}}}{\mbox{\texttt{min}}}$} $+$ \Ovalbox{\color{blue}$6 \frac{\mbox{ \texttt{cherry}}}{\mbox{\texttt{min}}}$}

\textbf{Question 4.} In what time they will fill the bowl?

\textbf{Answer:}  \Ovalbox{\color{blue}$6 \;\mbox{\texttt{min}}$} = \Ovalbox{\color{blue}$48 \mbox{ \texttt{cherry}}$} $\div$ \Ovalbox{\color{blue}$8\frac{\mbox{ \texttt{cherry}}}{\mbox{\texttt{min}}}$}

\end{boxedminipage} \\\begin{boxedminipage}[t]{12cm}
\

\begin{center}
 The answer to the problem is\\
Working together, Alice and Bob will fill the bowl in {\color{green} \Ovalbox{$6\; \mbox{\texttt{min}}$}}\\

\end{center}

\end{boxedminipage} \\

\begin{boxedminipage}[t]{12cm}
\
You may click on highlighted red {\color{red}\Ovalbox{numbers}} and change them! You will get a different problem of the same type.
\hfill\hamburger

\end{boxedminipage} \\
\end{center}
\caption{Manipulation with `helpful number': $72 \mbox{ \texttt{cherry}}$ replaced by $48 \mbox{ \texttt{cherry}}$---the answer did not change.}
\label{fig:screen-3}
\normalsize
\end{center}
\end{figure}

\begin{figure}
\begin{center}
\footnotesize
\bigskip
\begin{center}
\begin{boxedminipage}[t]{12cm}
\begin{center}
{\color{blue}{\large\textsc{Bowl of cherries}}}
\end{center}
\end{boxedminipage}
\end{center}
\begin{center}
\begin{boxedminipage}[t]{12cm}
\

\textbf{Data:} Bob fills the bowl in {\color{red} \Ovalbox{$8\; \mbox{\texttt{min}}$}}, Alice in {\color{red} \Ovalbox{$A\; \mbox{\texttt{min}}$}}.\\
$\mbox{ }\quad$ $\mbox{ }$   Now they work together.

\textbf{Question:}  In what time they will fill the bowl?
\end{boxedminipage}
\begin{boxedminipage}[t]{12cm}

\textbf{Helpful number: }The bowl contains {\color{red} \Ovalbox{$1\; \mbox{\texttt{cherry}}$}}.

\end{boxedminipage} \\

\begin{boxedminipage}[t]{12cm}
\

\noindent
\textbf{Question 1.} What is Alice's picking speeds?

\textbf{Answer:} \Ovalbox{\color{blue}$\frac{1}{A} \frac{\mbox{ \texttt{cherry}}}{\mbox{\texttt{min}}}$} = \Ovalbox{\color{red} $1 \mbox{ \texttt{cherry}}$} $\div$ \Ovalbox{\color{red}$A\; \mbox{\texttt{min}}$}

\textbf{Question 2.} What  is Bob's picking speed?

\textbf{Answer:}  \Ovalbox{\color{blue}$\frac{1}{8} \frac{\mbox{ \texttt{cherry}}}{\mbox{\texttt{min}}}$} = \Ovalbox{\color{red} $1\mbox{ \texttt{cherry}}$} $\div$ \Ovalbox{\color{red}$8\; \mbox{\texttt{min}}$}

\textbf{Question 3}. What is their picking speed working together?

\textbf{Answer:}  \Ovalbox{\color{blue}$\frac{A+8}{8A} \frac{\mbox{ \texttt{cherry}}}{\mbox{\texttt{min}}}$} = \Ovalbox{\color{blue}$\frac{1}{A} \frac{\mbox{ \texttt{cherry}}}{\mbox{\texttt{min}}}$} $+$ \Ovalbox{\color{blue}$\frac{1}{8} \frac{\mbox{ \texttt{cherry}}}{\mbox{\texttt{min}}}$}

\textbf{Question 4.} In what time they will fill the bowl?

\textbf{Answer:}  \Ovalbox{\color{blue}$\frac{8A}{A+8} \;\mbox{\texttt{min}}$} = \Ovalbox{\color{blue}$1 \mbox{ \texttt{cherry}}$} $\div$ \Ovalbox{\color{blue}$\frac{A+8}{8A}\frac{\mbox{ \texttt{cherry}}}{\mbox{\texttt{min}}}$}

\end{boxedminipage} \\\begin{boxedminipage}[t]{12cm}
\

\begin{center}
 The answer to the problem is\\
Working together, Alice and Bob will fill the bowl in {\color{green} \Ovalbox{$\frac{8A}{A+8} \; \mbox{\texttt{min}}$}}\\

\end{center}

\end{boxedminipage} \\

\begin{boxedminipage}[t]{12cm}
\
You may click on highlighted red {\color{red}\Ovalbox{numbers}} and change them! You will get a different problem of the same type.
\hfill\hamburger

\end{boxedminipage} \\
\end{center}
\caption{Substitution of a letter for a datum.}
\label{fig:screen-4}
\normalsize
\end{center}
\end{figure}

\section{Named numbers} \label{sec:Named}

\subsection{Types of named numbers}
Merging Algebra with Computer Science requires development of an appropriated  Domain Specific Language which should support calculations with \emph{named numbers}, of the kind as in the situation when we \emph{share} 10 apples among 5 people:
\[
10\, \mbox{\texttt{apples}}\, \div 5\, \mbox {\texttt{people}} = 2 \frac{\mbox{\texttt{apples}}}{\mbox{\texttt{people}}}.
\]
Obviously, we need a typed language where variables could be of arbitrary type assigned to them, say \texttt{apple}, or \texttt{people}, where only  variables of commensurable (or comparable) types could be added, say $1 \mbox {\texttt{m}} + 1 \mbox{ \texttt{cm}} = 101 \mbox{ \texttt{cm}}$,  but not apples and people; this restriction is possible to achieve\footnote{See Sirotin \cite{Sirotin2022}, a library of functions and objects of programming language  \textsc{Kotlin} (\url{https://kotlinlang.org/}) which allows working with variables whose values are expressed in the International System of Units \textsc{SI}  like meter, second etc. as well as some other common units like currencies, percentages etc..}; however, division of \texttt{apples} by \texttt{people} automatically produces a variable of a different type, $\frac{\mbox{\texttt{apples}}}{\mbox{\texttt{people}}}$.  Some types, such as \textsc{S\i} units or currencies,  need to be standard, built in the DSL---but Learner has to be able to create any types of fancy, say \texttt{ducat}  or \texttt{piastre} and declare exchange rate $1 \mbox{ \texttt{ducat}} == 5 \mbox{ \texttt{piastre}}$, so the types \texttt{ducat} and \texttt{piastre} become commensurable.

A more sophisticated manipulation with types is needed when we \emph{dispense} $10$ apples, $2$ apples a person, and wish to know how many people will get their apples:
\[10\, \mbox{\texttt{apples}}\, \div
\, 2\, \frac{\mbox{\texttt{\texttt{apples}}}}{\mbox{\texttt{\texttt{people}}}} = 5 \, \mbox{\texttt{\texttt{people}}}.
\]

\subsection{A bit of history}
Actually, it was clearly understood by the father of modern (symbolic) algebra, Fran\c{c}ois Vi\'{e}te who in 1591  wrote in his
\emph{Introduction to the Analytic Art} \cite[p. 16]{Vieta}  that

\small\bq\noindent
If one magnitude is divided by another, [the quotient] is
heterogeneous to the former \dots\  Much of the fogginess
and obscurity of the old analysts is due to their not paying
attention to these [rules].\eq\normalsize

Alas, these words remain true in the XXI Century.

This is the dirty secret of school arithmetic: it lives in an intimate relationship with the arithmetic of \emph{types} which is carefully hidden from children (and from many teachers). However, this is well known in physics where types are called \emph{dimensions}, and where \emph{dimensional analysis} is a simple, but powerful method of understanding of relations between quantities and magnitudes of different nature, which, in particular, gives a way to produce quick, frequently even quick \emph{mental} estimates of magnitudes. For example, \cite[Section 8.4]{Borovik-MuM}  contains a one page deduction of the legendary Kolmogorov's ``5/3'' law for the  energy spectrum of turbulent movement of gas or liquid.

\subsection{An elementary example from physics}
Let us look at one of simpler applications of dimension (type) analysis in physics.

Galileo Galilei observed that the period of pendulum does not depend on the  amplitude (span) of its movements.

So, the period $T$ (of type \texttt{second}) of pendulum depends on its length $L$ (of type \texttt{meter}) and acceleration of gravity $g$ (of type $\frac{\mbox{\texttt{meter}}}{\mbox{\texttt{second}}^2}$).

The only type-consistent formula which can be made out of that is
\[
T = C\sqrt{\frac{L}{g}},
\]
where $C$ is dimensionless constant. Even if we do not know the value of $C$ (this requires more subtle arguments), we may draw very interesting consequences.

For example, bipodal walking can be modeled as a sequence of falls in which a leg (of length $L$) which moves forward behaves as a pendulum of period $T$. Hence the speed $V$ is proportional to $\frac{L}{T}$ and
\[
V \sim \frac{L}{T} \sim \frac{L}{\sqrt{\frac{L}{g}}} \sim \sqrt{gL}.
\]

Hence walking on stilts increases speed (anyone who tried to walk on stilts  knows that),   while walking on the Moon is $\sqrt{6}$ slower than on Earth (since the acceleration of gravity on the Moon is about $\frac{1}{6}g$).

Dimensional analysis should be part of the school course of Physics. It is simple, beautiful, and has a potential to produce revelations.

\subsection{Type analysis in arithmetic}

Let us apply the type analysis to a classical old problem,  a part of mathematical folklore:
\bq
\textsc{Rabbits and Chicken.}  Mary has pets, some rabbits and some chicken. Together her pets have
$12$ heads and $32$ legs.
\

How many rabbits and how many  chickens does Mary have?
\eq

It is more open to type analysis because it obviously involves more than one type.

First of all, we have to carefully assign types to every piece (datum) of quantitative data that we are given, or know from out lived experience.

Of course, we are given data of  types \texttt{leg} and  \texttt{head}. So, we have to introduce variables \textsc{Legs} and \textsc{Heads}  of these types, respectively, and we are given their values.

Since every rabbit and every chicken has exactly one head, we have variables \textsc{Rabbits} and \textsc{Chicken} which could be safely thought of being of type  \texttt{head}.

Not much so far. Perhaps we have to turn to our lived experience: comparative anatomy of rabbits and chicken. Perhaps it means
\[
\mbox{\textsc{RabbitAnatomy}} = 4 \frac{\mbox{\texttt{leg}}}{\mbox{ \texttt{head}}}
\]
and
\[
\mbox{\textsc{ChickenAnatomy}} = 2 \frac{\mbox{\texttt{leg}}}{\mbox{ \texttt{head}}}
\]
So all data that we have are
\bea
 \mbox{\textsc{Heads}} &=& 12 \,\mbox{\texttt{head}}\\
 \mbox{\textsc{Legs}} &=& 32\, \mbox{\texttt{leg}}\\
\mbox{\textsc{RabbitAnatomy}} &=& 4 \frac{\mbox{\texttt{leg}}}{\mbox{ \texttt{head}}}\\
\mbox{\textsc{ChickenAnatomy}} &=& 2 \frac{\mbox{\texttt{leg}}}{\mbox{ \texttt{head}}}
\eea
We have to understand, that for a child, this analysis of data is challenging -- but exactly this skill,
\bq
\emph{ability to see mathematical structures and relations in the real world}
\eq
-- is missing in the mainstream mathematics education. But it has to be taught, and systematically.

Taking into account Fran\c{c}ois Vi\'{e}te's advice, and the  natural limitations of type arithmetic, we have \emph{very little choice} of sensible arithmetic operations between these values, we emphasize, \emph{we have very little choice}. One of the \emph{very few} actions that we may try is to compute
\[
 \mbox{\textsc{Heads}} \times \mbox{\textsc{RabbitAnatomy}} = 12 \,\mbox{\texttt{head}} \times 4 \frac{\mbox{\texttt{leg}}}{\mbox{ \texttt{head}}} = 48\, \mbox{ \texttt{leg}}.
\]
What we see here is \emph{type inference}, in terminology of computer science. But does it have a real life meaning? Yes, it does. This is the number of legs Mary's pets would have it all of them were rabbits.  This is more  than the given number of legs. By how much more?
\[
 \mbox{\textsc{ExcessiveLegs}} =  \mbox{\textsc{Heads}} \times \mbox{\textsc{RabbitAnatomy}}  -  \mbox{\textsc{Legs}} =  16\, \mbox{\texttt{leg}}.
\]
Where do these excessive legs come from? From chicken, each getting extra
\[
\mbox{\textsc{RabbitAnatomy}} - \mbox{\textsc{ChickenAnatomy}}  = 4 \frac{\mbox{\texttt{leg}}}{\mbox{ \texttt{head}}} - 2 \frac{\mbox{\texttt{leg}}}{\mbox{ \texttt{head}}} = 2 \frac{\mbox{\texttt{leg}}}{\mbox{ \texttt{head}}}
\]
legs. So, what is the number of chicken?
\bea
\mbox{\textsc{Chicken}} &=&  \frac{\mbox{\textsc{ExcessiveLegs}}}{\mbox{\textsc{RabbitAnatomy}} - \mbox{\textsc{ChickenAnatomy}}}\\[1ex] &=&
16\, \mbox{\texttt{leg}} \div  2 \frac{\mbox{\texttt{leg}}}{\mbox{ \texttt{head}}}\\[1ex]  &=& 8 \mbox{ \texttt{head}}
\eea
Of course, the number of rabbits after that is obvious:
\[
\mbox{\textsc{Rabbits}} = \mbox{\textsc{Heads}} - \mbox{\textsc{Chicken}} = 12 \,\mbox{\texttt{head}} - 8\,\mbox{\texttt{head}} = 4\,\mbox{\texttt{head}}.
\]
One of the \emph{very few} other available options to start our solution could be
\[
 \mbox{\textsc{Heads}} \times \mbox{\textsc{ChickenAnatomy}} = 12 \,\mbox{\texttt{head}} \times 2\frac{\mbox{\texttt{leg}}}{\mbox{ \texttt{head}}} = 24\, \mbox{ \texttt{leg}},
\]
but it leads to essentially the same solution.

This is not the shortest way to solve this problem but it has an advantage of being very formal, with every step forced on us, and leading to a computer code manipulating with named numbers -- this solution is already almost a code.

Alternative approaches are possible. Recently, the first author had one of his regular chats on \textsc{Google Meet} with a small group of Ukrainian refugee children (of ages about 8 to 11), and discussed with them this problem, and one of the kids suggested to look at \texttt{extremity}, or a limb, that is,  \texttt{leg} or \texttt{wing}, which leads to a more efficient solution -- but only because a deeper insight in the real world is applied, and more types are used.  Indeed, just three straightforward questions:
\begin{enumerate}
  \item What is the number of extremities?
  \item What is the number of wings?
  \item What is the number of chicken?
\end{enumerate}
lead to a solution.



\section{We need fusion of algebra and computing, not concatenation of the two}

A paper by Ian Benson and Jim Thorpe \cite{Benson-Thorpe2022} advertises the old paper by Trevor Fletcher \cite{Fletcher1971} \emph{Thinking with Arrows} as `conceptual mathematics', and a way to `represent and reason' about mathematical patterns.

Perhaps driven by the law of nominative determinism, Fletcher was enthusiastic about use of arrows in school mathematics. This is a fragment from his paper given by Benson and Thorpe as an example of `thinking with arrows'.

\begin{quotation}
If you have a tap that fills a bath in $3$ minutes and a tap that fills the same bath in $6$ minutes, you have to know how to decide that together they would take $2$ minutes. Likewise if you are told $10$ minute and $10$ minutes you have to be able to reply $5$ minutes. How is the calculation done? This time the operation of getting $2$ from $3$ and $6$ may be called \emph{tap}, and the associated diagram is

\begin{equation}\label{eq1}
\begin{diagram}
(3,6) &&\rTo^{tap} && 2\\
\dDotsto^{recip} && &&\dDotsto_{recip}\\
&&&& \\
\left(\frac{1}{3}, \frac{1}{6}\right) &&\rTo^{+} && \frac{1}{2}
\end{diagram}
\end{equation}
\end{quotation}

Here, we do not see how the answer to the problem was \textbf{found}; we see only that the \textbf{answer} to the problem is \textbf{represented} by a commutative diagram which reduces the answer (treated as a binary operation) to the binary operation of addition of positive real numbers; please notice that the word `operation' is used by Fletcher.

We have  objections to this approach.
\bi
\item We think the aim of mathematics education should be teaching the child the art of \textbf{solving} problems and \textbf{getting answers}.

\item \textbf{Interpretation} of answers should play important role; however, insistence on a particular form of \textbf{representation} of the answer is undesirable; after all

\item There should be a smooth transition between the arithmetic / algebraic thinking in solving a problem and computer science / computer programming thinking in representing the solution as an algorithm.

\item We need fusion of algebra and computing, not concatenation---but Diagram \ref{eq1} is a concatenation of two disjoint conceptual areas.
\ei

We wish to demonstrate, that, in the case of a class of problems on the arithmetic / algebra boundary where the tap problem naturally belongs, the commutative diagram representation of a solution is misleading and counterproductive.

Indeed, consider the `tap problem' together with three other problems from the same stable: proportionality dependencies between time, speed, and distance (where `distance' could perhaps be `level of water in the bath').

\bd
\item[Problem 1]  If you have a tap that fils a bath in $a$ minutes and a tap that fills the same bath in $b$ minutes, in how many minutes
 they  would fill the bath together? This is, of course, Fletcher's problems with answer represented  in Diagram 1.
 \item[Problem 2]  If a car travels from $A$ to $B$ at speed of  $a$ miles per hour, and then returns from $B$ to $A$ at speed of $b$ miles per hour, what is the average speed of the car on the whole journey?
 \item[Problem 3] Two cars started their journeys at sunrise, one from $A$ to $B$, and another one from $B$ to $A$. They met at noon and completed their journeys at $a$  hours p.m.\ and $b$ hours p.m., correspondingly. By how much was sunrise earlier than the noon on that day?
 \item[Problem 4]  There are two cities on a river, one is upstream of another. It takes for a steamboat $a$ days to get from one city to another, and $b$ days to return back. How many days it will take a raft to drift from the city which is upstream to the downstream city?
\ed
We give here answers to these problems in random order:
\begin{center}
(a) $\sqrt{ab}$;\; (b) $\displaystyle \frac{ab}{a+b}$;\; (c) $\displaystyle \frac{2ab}{|a-b|}$;\; (d) $\displaystyle \frac{2ab}{a+b}$;
\end{center}
How would Learner know which answer is for which problem? Being told by the teacher? Where do these answers come form? Fletcher's commutative diagram is indeed commutative, but, we are afraid, this is a bad example of theoretical-methodology-driven approach to teaching, not a Learner-focused approach to Learner's learning. On the contrary, we suggest, for the critical stage of moving from arithmetic to algebra, a smooth way from an algorithmic solution of an arithmetic problem to its algebraic expression and coding. Commutative diagrams could come into the play a few years later, not now.

If we go into the underlying mathematics, we quickly discover that only one of these answers (a)-(d) can be represented by a commutative diagram similar to Diagram \ref{eq1}. Indeed, if $\circ$ and $\diamond$ are two binary operations on the set of positive real numbers $\mathbb{R}^{>0}$ satisfying the commutative diagram
\begin{equation}\label{eq2}
\begin{diagram}
(a,b) &&\rTo^{\circ} && a\circ b\\
\dDotsto^{\phi} && &&\dDotsto_{\phi}\\
&&&& \\
\left(\phi(a), \phi(b) \right) &&\rTo^{\diamond} && \phi(a) \diamond \phi(b)
\end{diagram}
\end{equation}
for some map $\phi: \mathbb{R}^{>0}\longrightarrow  \mathbb{R}^{>0}$, then
\begin{equation}\label{eq3}
\phi(a\circ b) = \phi(a)  \diamond \phi(b),
\end{equation}
that is, $\phi$ is a homomorphism from monoid $A= \{\mathbb{R}^{>0}, \circ\}$ to monoid $B=\{\mathbb{R}^{>0}, \diamond\}$. In answers (a) and  (d)---by the way, these are the geometric and harmonic means, quite important concepts on their own---the monoid $A$ is \emph{idempotent}, that is, $a \circ a = a$ for all $a\in A$. But then the image  $\Im(\phi)$ of $\phi$ is also idempotent. Hence $\diamond$ cannot be, for example, the standard addition because $\{\mathbb{R}^{>0}, +\}$ has no idempotents,   and if  $\diamond$ is the standard multiplication then $\Im(\phi)=\{1\}$ and $\phi(a) = 1$ for all $a\in A$. This is the trivial homomorphism which carries no information about the operation $\circ$.

Answer (c) is not even an operation: it is not defined at $a=b$. It has a singularity there, which makes it even more difficult to represent it by a commutative diagram of kind of Diagram  \ref{eq1} with a nice, and familiar to Learner, binary operation $\diamond$ in the arrow at the bottom.

Equation \ref{eq3} provides a way of constructing an unbelievable variety of commutative diagrams of the kind of Diagram \ref{eq2}. You have to start with an arbitrary binary operation $\diamond$  on  $\mathbb{R}^{>0}$ and take an arbitrary 1--1 bijection $\phi$ on (permutation of) $\mathbb{R}^{>0}$ and define
\begin{equation}\label{eq4}
 a\circ b = \phi^{-1}(\phi(a) \diamond \phi(b)),
\end{equation}
then Diagram \ref{eq2} is obviously commutative. And we are spoilt for choice: there are at least $2^{2^{\aleph_0}}$ permutations $\phi$. The only question is
\bq
How many of these commutative diagrams are relevant to school arithmetic and elementary school algebra?
\eq
We are afraid, just a handful. Some of them are interesting, but still rather isolated.

For example, we can take $\diamond = +$, and $\phi(x) = x^n$, $n$ a natural number.
We get a commutative diagram
\begin{equation}\label{eq5}
\begin{diagram}
(a,b) &&\rTo^{\circ} && \sqrt[n]{a^n +  b^n}\\
\dDotsto^{x\mapsto x^n} && &&\dDotsto_{x\mapsto x^n}\\
&&&& \\
\left(a^n, b^n \right) &&\rTo^{+} && a^n +b^n
\end{diagram}
\end{equation}
The case $n=2$ is exhibited by Fletcher as `Pythagorean map':
\begin{equation}\label{eq6}
\begin{diagram}
(3,4) &&\rTo^{\mbox{\texttt{Pythag}}} && 5\\
\dDotsto^{\mbox{\texttt{Sq}}} && &&\dDotsto_{\mbox{Sq}}\\
&&&& \\
\left(9, 16 \right) &&\rTo^{+} && 25
\end{diagram}
\end{equation}
---but the diagram gives no indication of how the Pythagoras' theorem has been proved.

\section{Conclusion}

In this brief discussion, it has to be emphasised:
\medskip

\begin{center}
\shadowbox{\parbox{10cm}{%
Questions in the solution are supposed to be formulated by \textbf{Learner}, not by a teacher. \textbf{Learner has to be in control.} \cite{Borovik2017-In-control}
}}
\end{center}
\medskip

What we have seen:

\bi
\item Careful analysis of data.
\item A sequence of `questions' which naturally becomes an algorithm.
\ei
Also, concepts from Computer Science:
\bi
\item Typed variables, aka `named numbers'.
\item Evaluation by call by name and by call by value.
\item Reification as a way of introduction of intermediate parameters.
\ei

\section{Some socio-political comments}

\begin{flushright}
  \emph{L'enfer, c'est les autres.}\\
  Jean-Paul Sartre, \emph{Huis Clos}, 1943
\end{flushright}

The old (1996!) paper by Neil Koblitz\footnote{Neil Koblitz is one of the founding fathers of algebraic cryptography, he proposed the elliptic curve cryptography.} \emph{The case against computers in K-13 math education} \cite{Koblitz1996} is well worth re-reading. His warnings are still relevant and fully apply to the possibility of badly prepared, underfunded, hastily introduced computerisation of mathematics education in schools:
\emph{\bq
The downside can be divided into several broad areas:
\bi
\item drain on resources (money, time, energy);
\item bad pedagogy;
\item anti-intellectual appeal;
\item corruption of educators.
\ei
\eq}
We hope this paper may serve as a warning: proper introduction of  ``computational thinking'' in school mathematics is an immensely complex and highly expensive task.

We can only applaud the heroic efforts of hundreds of teachers in various countries who are trying to do something in this direction, and we think it is very important to help them in every way available, their experience needs to be studied,  systematised, and fed into the development of a potential reform. Very often they are some of the best teachers in their countries. Unfortunately, this usually means that methods and approaches discovered and developed by them are, by default, not scalable and not transferable to the entire country, and for a basic reason: most other teachers are not like them, they are less educated, less motivated, and frequently demoralised by the daily grind of school work.

A serious reform requires a systematic re-education of a whole army of teachers and providing them with (paid, of course) time for personal professional development. An immediate consequence: we will need more teachers.

Lessons of Kolmogorov's reform of school mathematics curriculum in Russia in the 1970s \cite{Borovik2022} need to be learnt. Andrey Kolmogorov was the world famous mathematician, he had best intentions,  but his reform spectacularly failed, and one of the reasons for that was underestimating the role of teachers, even direct neglect of teachers. What is really sad is that Russia at that time had a well developed and functional system of preparing mathematics teachers via a numerous pedagogical colleges (4 years or tertiary education) and mathematics departments in regional universities (5 years)---but this resource had not been engaged. Arguably, this educational powerhouse degenerated over the years, but still the situation in Russia is unlikely to be as bad as that in the UK where, in the words of Tony Gardiner (the best expert in mathematics education in Britain) \cite{Gardiner2018},
\bq
\emph{We know of no other system that pretends to train mathematics teachers by placing small groups of trainees at the mercy of teachers with no
relevant ITE [Initial Teachers' Education] experience beyond being themselves teachers, with much of the input being ``generic'' rather than subject-specific. England appears to be alone among developed nations in embracing such an approach.}
\eq

Recalling again the tragic fate of Andrei Kolmogorov and his reform \cite{Borovik2022}, we reiterate our warning:

\medskip

\begin{center}
\shadowbox{\parbox{10cm}{%
Any attempt of a deep reform of mainstream mathematics education is a huge task. It is dangerous to undertake it lightly.
}}
\end{center}

\section*{Acknowledgements}

The first author is forever grateful to Martin Hyland who encouraged him to get involved in this project, and  thanks Ian Benson, Tony Gardiner, and Victor Sirotin for many useful discussions and corrections. His thanks go to Rick Booth, Anja Meyer, and Natasa Strabic for their understanding, support and helpful comments.

\section*{Disclaimer} The authors write in their personal capacities and the views expressed do not necessarily represent the position of their employers or any other person, corporation, organisation, or institution.\\

\section*{References}
\begingroup
\renewcommand{\section}[2]{}%

\endgroup

\bigskip


\end{document}